# КУБИЧЕСКИЕ СИСТЕМЫ ТИПА ДАРБУ С НЕЭЛЕМЕНТАРНОЙ ОСОБОЙ ТОЧКОЙ НА ЭКВАТОРЕ ПУАНКАРЕ

Е. П. ВОЛОКИТИН

**Аннотация.** Мы исследуем глобальное поведение траекторий полиномиальной системы $\dot x = x - x^2y + pxy^2 + y^3$, $\dot y = y + py^3$, $p \in \mathbb{R}$. Наше исследование примыкает к работе [6].

**Ключевые слова:** полиномиальные системы, особые точки, экватор Пуанкаре, фазовые портреты, рациональные интегралы.

## Введение

Рассмотрим плоскую систему обыкновенных дифференциальных уравнений

(1) $$\dot x = x + P_n(x,y), \ \dot y = y + Q_n(x,y),$$

где $P_n(x,y)$, $Q_n(x,y)$ — однородные многочлены $n$-й степени. Мы будем называть такую систему системой типа Дарбу.

Системы вида (1) рассматривались различными авторами, см. [1]–[6] и процитированную там литературу. В этих работах изучались традиционные для качественной теории дифференциальных уравнений вопросы такие, как интегрируемость, наличие или отсутствие предельных циклов, локальные и глобальные фазовые портреты и т. д.

В [6] исследовались общие свойства систем вида (1), и в качестве примера применения полученных результатов были рассмотрены кубические системы типа Дарбу. Для них авторы получили полный по их мнению список глобальных фазовых портретов (с точностью до топологической эквивалентности). Однако вопреки утверждению авторов этот список не является полным. Для оправдания нашего утверждения мы рассмотрим однопараметрическое семейство кубических систем типа Дарбу вида

(2) $$\dot x = x - x^2y + pxy^2 + y^3, \ \dot y = y + py^3, \ p \in \mathbb{R}$$

и предъявим фазовые портреты систем этого семейства, не содержащиеся в списке, полученном в [6].

Одной из причин допущенной авторами [6] ошибки является тот факт, что они не учли возможности существования у систем вида (1) неэлементарной особой точки на экваторе Пуанкаре с равной нулю матрицей линейного приближения ("линейный ноль"). В этом случае возможно строение окрестности этой точки, которое имеет более сложную структуру по сравнению с традиционно возникающими и учтёнными в [6] негиперболическими точками типа узел и седло-узел. Мы рассматривали сходные ситуации в [5], [7]. В настоящей работе мы проведём соответствующее исследование с помощью более прозрачного и простого метода, применяя алгебраическое раздутие без использования трансцендентных уравнений. Кроме того мы найдём интеграл системы (2) и обсудим его свойства.

## Основная часть

Рассмотрим систему

(2) $$\dot{x} = x - x^2 y + pxy^2 + y^3, \ \dot{y} = y + py^3, \ \ p \in \mathbb{R}.$$

**Теорема 1.** *1) Фазовые портреты системы* (2) *приведены на рис. 3.*
*2) Если $p \neq 0$, система* (2) *имеет интеграл*

$$H(x,y) = (x-y)^p (x+y)^{-p}(1+py^2),$$

*который является рациональным, если $p \in \mathbb{Q} \setminus \{0\}$.*
*3) Если $p = 0$, система* (2) *имеет интеграл*

$$H(x,y) = \frac{x-y}{x+y} e^{y^2}.$$

*Если $p = 0$, система* (2) *не имеет рационального интеграла.*

*Доказательство.* Для изучения системы (2) мы применим более прозрачные и простые по сравнению с [5], [7] методы исследования. Для разрешения особенностей здесь мы будем использовать алгебраическое раздутие, что позволяет не использовать трансцендентные функции. Кроме того это даёт возможность оценить асимптотику траекторий. К тому же в дополнение к [5] мы найдём интегралы системы (2) и укажем их интересные на наш взгляд свойства.

Система (1) детально исследовалась в [1], [4], [6]. В частности, там получено, что поведение в целом траекторий системы (2) определяется поведением её траекторий вблизи экватора Пуанкаре.

Экватор Пуанкаре отвечает оси $\{z=0\}$ системы

(3) $$\dot{u} = Q_3(1,u) - uP_3(1,u), \ \dot{z} = -z(z^2 + P_3(1,u)),$$

которая является результатом компактификации Пуанкаре системы (2) и получается из неё с использованием замены

$$u = \frac{y}{x}, \ z = \frac{1}{x}.$$

Подробнее см. [8], [9].

В нашем случае компактификация системы (2) даёт систему

(4) $$\dot{u} = u^2(1-u^2), \ \dot{z} = -z(z^2 - u + pu^2 + u^3),$$

Система (4) имеет на экваторе три особые точки $O_0(0,0), O_1(1,0), \ O_2(-1,0)$. Матрицы линейного приближения в этих точках суть

(5) $$J_0 = \begin{pmatrix} 0 & 0 \\ 0 & 0 \end{pmatrix}, \ J_1 = \begin{pmatrix} -2 & 0 \\ 0 & -p \end{pmatrix}, \ J_2 = \begin{pmatrix} 2 & 0 \\ 0 & -p \end{pmatrix}.$$

Мы видим, что начало координат $u = 0, z = 0$ является неэлементарной особой точкой: матрица линейной части нулевая. Такие точки были упущены из рассмотрения в [6].

Для изучения особой точки $O_0(0,0)$ сделаем замену координат (раздутие) [10], [11]

(6) $$u = u, w = \begin{cases} \frac{z}{\sqrt{u}}, & u > 0, \\ \frac{z}{\sqrt{-u}}, & u < 0. \end{cases}$$

После применения раздутия получим

(7) $$\dot{u} = u^2(1-u^2), \ \dot{w} = \frac{1}{2}uw - uw^3 - pu^2w - \frac{1}{2}u^3w, \quad u > 0,$$

(8) $$\dot{u} = u^2(1-u^2), \ \dot{w} = \frac{1}{2}uw + uw^3 - pu^2w - \frac{1}{2}u^3w, \quad u < 0.$$

Введём новое время $d\tau = udt$. Получим системы с теми же траекториями

(9) $$\dot{u} = u(1-u^2), \ \dot{w} = \frac{1}{2}w - w^3 - puw - \frac{1}{2}u^2w, \quad u > 0,$$

(10) $$\dot{u} = u(1-u^2), \ \dot{w} = \frac{1}{2}w + w^3 - puw - \frac{1}{2}u^2w, \quad u < 0.$$

Дифференцирование по новой переменной по-прежнему обозначено точкой.

Система (9) имеет на оси $\{u = 0\}$ три гиперболические особые точки $\bar{O}_0(0,0), \bar{O}_1(0, 1/\sqrt{2}), \bar{O}_2(0, -1/\sqrt{2})$, у которых матрицы линейной части имеют вид

$$\bar{J}_0 = \begin{pmatrix} 1 & 0 \\ 0 & 1/2 \end{pmatrix}, \ \bar{J}_1 = \begin{pmatrix} 1 & 0 \\ -p/\sqrt{2} & -1 \end{pmatrix}, \ \bar{J}_2 = \begin{pmatrix} 1 & 0 \\ p/\sqrt{2} & -1 \end{pmatrix}.$$

Точка $\bar{O}_0(0,0)$ будет неустойчивым узлом, оставшиеся две — сёдлами.

Система (10) имеет на оси $\{u = 0\}$ единственную особую точку — гиперболический устойчивый узел.

На рис.1 изображены траектории систем (7), (8) в окрестности оси $\{u = 0\}$.

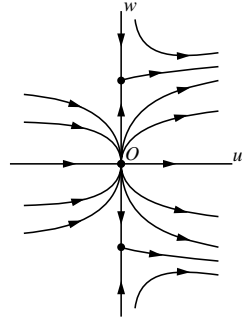

Рис. 1. Траектории систем (7), (8) в окрестности оси $\{u = 0\}$

Обратная замена координат (сжатие) позволяет изобразить поведение траекторий системы (4) в окрестности неэлементарной особой точки $O_0(0,0)$. Эта окрестность содержит два узловых и два седловых сектора, рис. 2.

Особые точки $O_1(1,0), \ O_2(-1,0)$ системы (4) имеют следующий характер (см. (5)):

$p > 0, O_1(1,0)$ — гиперболический устойчивый узел,
$O_2(-1,0)$ — гиперболическое седло;
$p < 0, O_1(1,0)$ — гиперболическое седло,
$O_2(-1,0)$ — гиперболический неустойчивый узел.

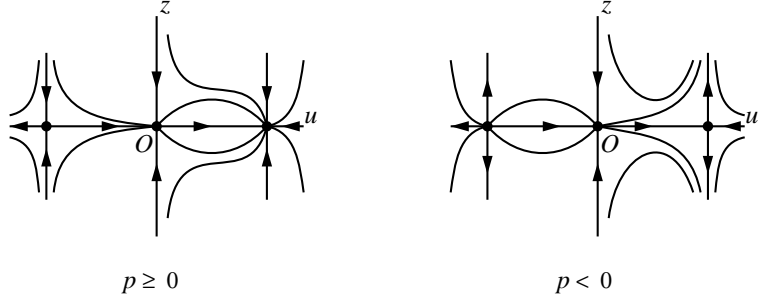

Рис. 2. Бесконечно удалённые особые точки системы (4)

В случае $p = 0$ точки $O_1(1,0), O_2(-1,0)$ будут элементарными полугиперболическими точками. Мы применяем разработанные для таких точек методы исследования, предложенные в [8], [9], и получаем, что $O_1(1,0)$ будет устойчивым узлом, $O_2(-1,0)$ — седлом.

На рис. 2 изображены варианты поведения траекторий системы (4) в окрестности оси $\{z = 0\}$, отвечающей экватору Пуанкаре, в зависимости от значений параметра $p$. Системы с таким поведением траекторий в окрестности экватора не попали в поле зрения авторов [6].

Охарактеризуем поведение траекторий системы (2) в конечной части плоскости. Мы будем опираться на результаты [4], [6].

Система имеет три инвариантные кривые $y = 0, y = x, y = -x$.

Напомним, что стационары системы (2) могут располагаться только на инвариантных прямых.

Вдоль прямой $y = 0$ (2) принимает вид

$$\dot{x} = x, \ \dot{y} = y,$$

и видно, что точки покоя, отличающиеся от начала координат, на этой прямой отсутствуют.

Вдоль прямой $y = x$ (2) принимает вид

$$\dot{x} = x + px^3, \ \dot{y} = x + px^3.$$

Если $p \geq 0$, на прямой отсутствуют точки покоя, отличающиеся от начала координат. Если $p < 0$, на прямой имеются две точки покоя $(\pm 1/\sqrt{-p}, \pm 1/\sqrt{-p})$, которые являются устойчивыми гиперболическими узлами.

Вдоль прямой $y = -x$ (2) принимает вид

$$\dot{x} = x + px^3, \ \dot{y} = -x - px^3.$$

Как и в предыдущем случае, если $p \geq 0$, на прямой отсутствуют точки покоя, отличающиеся от начала координат. Если $p < 0$, на прямой имеются две точки покоя $(\pm 1/\sqrt{-p}, \mp 1/\sqrt{-p})$, которые являются гиперболическими сёдлами.

Предельные циклы в системе (2) очевидно отсутствуют.

Суммируя полученную информацию, строим фазовые портреты системы (2) на диске Пуанкаре, рис. 3.

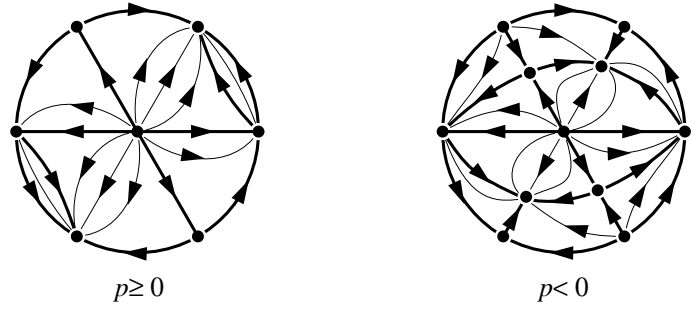

Рис. 3. Фазовые портреты системы (4)

Как уже отмечалось, такие фазовые портреты не могли быть обнаружены в процессе исследования [6].[1]

Пункт 1) теоремы доказан.

Используя метод Дарбу интегрирования систем ОДУ с помощью инвариантов, см. например, [12], мы можем найти интеграл (2).

Функция $L(x,y)$ называется инвариантом системы
$$\dot{x} = P(x,y), \dot{y} = Q(x,y),$$
если она удовлетворяет условию
$$DL \equiv \frac{\partial L(x,y)}{\partial x}P(x,y) + \frac{\partial L(x,y)}{\partial y}Q(x,y) = k(x,y)L(x,y),$$
где $k(x,y)$ — многочлен от переменных $x, y$, который называется кофактором инварианта $L(x,y)$.

Если система имеет инварианты $L_1, L_2, \ldots, L_s$ с кофакторами $k_1, k_2, \ldots, k_s$ и $\alpha_1 k_1 + \alpha_2 k_2 + \ldots \alpha_s k_s = 0$, функция $H = L_1^{\alpha_1} L_2^{\alpha_2} \ldots L_s^{\alpha_s}$ является интегралом системы.

Система (2) имеет инвариантами многочлены $L_1 = y$, $L_2 = x - y$, $L_3 = x + y$, $L_4 = 1 + py^2$ и экспоненциальный множитель $L_5 = e^{y^2}$ с кофакторами $k_1 = 1 + py^2$, $k_2 = 1 - xy - y^2 + py^2$, $k_3 = 1 - xy + y^2 + py^2$, $k_4 = 2py^2$, $k_5 = 2y^2(1 + py^2)$.

Имеем $pk_2 - pk_3 + k_4 = 0$.

В таком случае

(11) $$H(x,y) = L_2^p L_3^{-p} L_4 = (x-y)^p (x+y)^{-p}(1 + py^2), \quad p \neq 0$$

является первым интегралом системы (2).

Очевидно, что при рациональных значениях $p \neq 0$ интеграл (11) будет рациональным. Числа $h_1 = 0, h_2 = \infty$ будут его замечательными значениями.

Пункт 2) теоремы доказан.

Если $p = 0$, $H(x,y) \equiv 1$. При $p = 0$ имеем $k_2 - k_3 + k_5 = 0$, и

(12) $$H_0(x,y) = \frac{x-y}{x+y} e^{y^2}$$

---

[1]В [5] приведены также другие варианты глобальных фазовых портретов кубических систем типа Дарбу, которых нет в [6].

является первым интегралом системы (2), когда $p = 0$

(13) $$\dot{x} = x - x^2 y + y^3, \ \dot{y} = y.$$

Интеграл (12), очевидно, не является рациональным. Более того, можно показать, что система (13) вообще не имеет рационального интеграла.

Чтобы убедиться в этом, рассмотрим кривую, которая является линией уровня $H_0(x, y) = 1$

$$\Gamma = \{(x, y) \in \mathbb{R} : \frac{x - y}{x + y} e^{y^2} = 1, \ y > 0\}.$$

Кривая $\Gamma$ — сепаратриса бесконечно удалённой неэлементарной особой точки, лежащей на конце оси $Ox$, рис. 4.

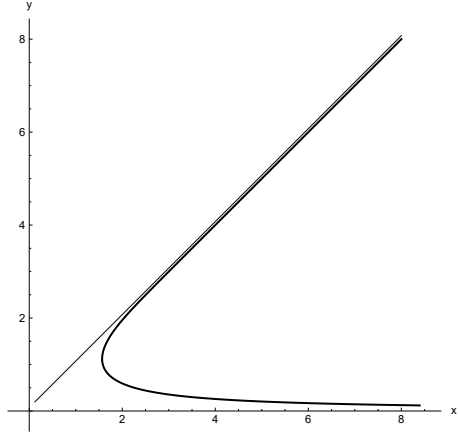

Рис. 4. Кривая $\Gamma = \{(x, y) \in \mathbb{R} : H_0(x, y) = 1, \ y > 0\}$

Пусть $\Gamma$ является ветвью алгебраической кривой $F(x, y) = 0$ степени $n$. Следуя Л. Эйлеру [13], запишем эту кривую в виде

(14) $$F(x, y) \equiv P(x, y) + Q(x, t) + R(x, y) + \cdots = 0,$$

где $P$ является высшим членом, содержащим в себе все члены степени $n$, $Q$ — вторым членом, содержащим члены степени $n - 1$, и так далее.

Кривая $F(x, y) = 0$ является инвариантной кривой системы (13). Согласно пункту 1) теоремы 1 неограниченные ветви могут иметь ассимптотами только прямые $y = 0, \ x - y = 0, \ x + y = 0$. В таком случае согласно Л. Эйлеру [13]

$$P(x, y) = y^k (x - y)^l (x + y)^m P_r(x, y), \ k + l + m + r = n,$$

где многочлен

$$P_r(x, y) \equiv \alpha x^r + \beta x^{r-1} y + \ldots \gamma y^r$$

не имеет линейных множителей. Инами словами уравнение $P_r(1, s) = 0$ не имеет действительных корней. В частности,

(15) $$P_r(1, -1) \equiv \alpha - \beta + \cdots + \gamma \neq 0.$$

Вдоль $\Gamma$ имеем

$$x = \frac{e^{y^2} + 1}{e^{y^2} - 1} y.$$

Обозначим $z = e^{y^2}$. Тогда вдоль Γ имеем

(16) $$x = \frac{z+1}{z-1}y, \quad \frac{x}{y} = \frac{z+1}{z-1}.$$

Подстановка (16) в (14) после очевидных преобразований приводит к равенству вида

(17) $$\mathcal{F}(y,z) = \mathcal{P}(y,z) + \mathcal{Q}(y,z) + \cdots = 0.$$

Здесь $\mathcal{F}(y,z)$, $\mathcal{P}(y,z)$, $\mathcal{Q}(y,z)$ — многочлены; при этом
$$P = y^k(x-y)^l(x+y)^m(\alpha x^r + \beta x^{r-1}y + \ldots \gamma y^r) =$$
$$y^n\left((\frac{x}{y}-1)^l(\frac{x}{y}+1)(\alpha(\frac{x}{y})^r + \beta(\frac{x}{y})^{r-1} + \cdots + \gamma)\right) =$$
$$y^n\left((\frac{2}{z-1})^l(\frac{2z}{z-1})^m(\alpha(\frac{z+1}{z-1})^r + \beta(\frac{z+1}{z-1})^{r-1} + \cdots + \gamma\right) =$$
$$2^{l+m}y^n z^m \frac{\alpha(z+1)^r + \beta(z+1)^{r-1}(z-1) + \cdots + \gamma(z-1)^r}{(z-1)^{l+m+r}},$$
$$\mathcal{P} = 2^{l+m}y^n z^m(\alpha(z+1)^r + \beta(z+1)^{r-1}(z-1) + \cdots + \gamma(z-1)^r)(z-1)^{n-l-m-r},$$
$$\mathcal{Q} = y^{n-1}q(z), \ldots.$$

Коэффициент $p_{nm}$ при мономе $y^n z^m$ в многочлене $\mathcal{P}$ имеет вид

(18) $$p_{nm} = (-1)^{n-l-m-r}2^{l+m}(\alpha - \beta + \ldots \gamma).$$

В силу (15) $p_{nm} \neq 0$, поэтому многочлен $\mathcal{F}$ не является тождественно нулевым.

В таком случае равенство (17) означает, что система функций $\{y^r e^{sy^2}\}_{r,s \in \mathbb{Z}}$ является линейно зависимой на промежутке $y > 0$, что очевидно неверно.

Итак, кривая Γ не алгебраическая. В таком случае система (13) не имеет рационального интеграла, поскольку при наличии рационального интеграла все траектории системы должны быть алгебраическими кривыми.

Пункт 3) теоремы доказан. $\square$

Отметим, что рациональный интеграл
$$H(x,y) = \frac{(x-y)(1+y^2)}{x+y},$$
полученный из (11) при $p = 1$, и интеграл
$$H(x,y) = \frac{(x-y)e^{y^2}}{x+y},$$
из (12), не являющийся рациональным, отвечают топологически эквивалентным фазовым портретам соответствующих систем.

## Список литературы

Волокитин Евгений Павлович
Институт математики им. С. Л. Соболева СО РАН,
пр. академика Коптюга, 4,
630090, Новосибирск, Россия.
*Email address*: volok@math.nsc.ru